\documentclass[12pt]{amsart}
\usepackage{amsmath, amssymb, latexsym, amsthm, mathrsfs, color, mathtools, soul}
\usepackage[hyphens]{url}
\usepackage[T1]{fontenc}
\usepackage{lmodern}

\usepackage[top=2in, bottom=1.5in, left=1in, right=1in]{geometry}

\newcommand{\nc}{\newcommand}

\nc{\thusfar}{\par\bigskip\centerline{\my{--- Edited thus far ---}}\par\bigskip}
\nc{\lei}{\le^\oo}
\nc{\card}[1]{\left|#1\right|}
\nc{\medcard}[1]{\biggl|\,#1\,\biggr|}
\nc{\smallcard}[1]{|\,#1\,|}
\nc{\bds}{bidirectional $\roth$-scale}
\nc{\bbN}{\mathbb{N}}
\nc{\beq}{\begin{eqnarray*}}\nc{\eeq}{\end{eqnarray*}}
\nc{\mbq}{\mb{?}}
\nc{\mb}[1]{{\mbox{\textbf{#1}}}}
\nc{\nop}{$\times$}
\nc{\fbn}{\!\!\fbox{\!\nop\!}\!\!}
\nc{\yup}{\checkmark}
\nc{\forces}{\Vdash}
\nc{\name}[1]{\dot{#1}}
\nc{\tf}{\my{FINISHED THUS FAR}}
\nc{\FU}{Fr\'echet--Urysohn}
\nc{\gs}{$\gamma$~space}
\nc{\Ga}{\Gamma}\nc{\Om}{\Omega}
\nc{\smallbinom}[2]{\begin{psmallmatrix} #1\\ #2 \end{psmallmatrix}}
\nc{\bgamma}{\smallbinom{\Om}{\Ga}}

\nc{\productive}[2]{\bigl(#1,\allowbreak #2\bigr)^\x}
\nc{\Sel}{\mathsf{S}}
\nc{\sset}[2]{\{\,#1 : #2\,\}}
\nc{\smb}[1]{{\!\!\mb{#1}\!\!}}
\nc{\medset}[2]{{\biggl\{\,#1 : #2\,\biggr\}}}
\nc{\smallmedset}[2]{{\bigl\{\,#1 : #2\,\bigr\}}}
\nc{\set}[2]{{\left\{\,#1 : #2\,\right\}}}
\nc{\seq}[2]{{\la\, #1 : #2\,\ra}}
\nc{\eseq}[1]{#1_0, \allowbreak #1_1, \allowbreak\dots} 
\nc{\cube}{(\Cantor)^\bbN}
\nc{\Match}{\op{Match}}
\nc{\concat}[1]{\hat{\phantom{a}}\langle #1\rangle}
\nc{\poset}{\mathbb{P}}
\nc{\fn}[1]{{\op{Fn}(#1\times\w,2)}}
\nc{\linadd}{\op{linadd}}
\nc{\nonprod}{\non^\x}
\nc{\alephes}{{\aleph_0}}
\nc{\my}[1]{\marginpar{\textcolor{red}{***}}\textcolor{red}{#1}}
\nc{\later}[1]{{\color{green} #1}}
\nc{\BTs}[1]{{\color{green} #1 (BT)}}
\nc{\Cp}{\op{C}_\mathrm{p}}
\nc{\Bp}{\op{B}_p}
\nc{\Pa}[8]{\bibitem{#1} {#2}, \emph{#3}, {#4} \textbf{#5} ({#6}), {#7}--{#8}.}
\nc{\tPa}[5]{\bibitem{#1} {#2}, \emph{#3}, {#4}, to appear.}
\nc{\sPa}[4]{\bibitem{#1} {#2}, \emph{#3}, {#4}, submitted.}
\nc{\Bc}[9]{\bibitem{#1} {#2}, \emph{#3}, in: \textbf{#4} (#5), #6 #7, #8--#9.}
\nc{\fD}{\mathfrak{D}}
\nc{\fX}{\mathfrak{X}}
\nc{\Onbd}{\Op_{\mathrm{nbd}}} 
\nc{\Omnb}{\Om_{\mathrm{nbd}}} 
\nc{\od}{\mathfrak{od}}
\nc{\Setting}[7]{\xymatrix@R=4pt@C=7pt{#1\ar@{-}[r]&#2\ar@{-}[r]&#3\\&#4\ar@{-}[u]\\
#5\ar@{-}[uu]\ar@{-}[r] & #6\ar@{-}[u]\ar@{-}[r] & #7\ar@{-}[uu]}}
\nc{\mx}[1]{\begin{matrix}#1\end{matrix}}
\nc{\plim}{p\txt{-}\lim}
\nc{\Bgp}{{\Z^\bbN}}
\nc{\Cgp}{{{\Z_2}^\bbN}}
\nc{\Cite}[1]{\textbf{[#1]}}
\nc{\Next}[1]{{#1^+}}
\nc{\cFin}{\mathrm{cF}}
\nc{\intvl}[2]{{[#1(#2),\allowbreak #1(#2\!+\!1))}}
\nc{\Bdd}{\mathbf{B}}
\nc{\Dfin}{\mathfrak{D}_\mathrm{fin}}
\nc{\grbl}{{\mbox{\textit{\tiny gp}}}}
\nc{\bbP}{\mathbb{P}}
\nc{\BOfat}{\B_{\Om_{\mathrm{fat}}}}
\nc{\Bgood}{\B_{\mathrm{good}}}
\nc{\compactN}{\cl{\mathbb{N}}}
\nc{\blocks}[2]{\op{cl}_{#2}(#1)}
\nc{\blocksplus}[2]{\op{cl}^+_{#2}(#1)}
\nc{\arx}[1]{\texttt{http://arxiv.org/math/#1}}
\nc{\bq}{\begin{quote}}
\nc{\eq}{\end{quote}}
\nc{\cl}[1]{\overline{#1}}
\nc{\CH}{the Continuum Hypothesis}
\nc{\MA}{Martin's Axiom}
\nc{\Bfat}{\B_\mathrm{fat}}
\nc{\inv}{^{-1}}
\nc{\Cantor}{{2^\w}}
\nc{\bP}{\mathbf{P}}
\nc{\bof}{\op{\fb}}
\nc{\dof}{\op{\fd}}
\nc{\bofF}{\bof(\cF)}
\nc{\sr}[3]{\underset{\mbox{#3}}{\mbox{#1}}}
\nc{\gp}{\binom{\Om}{\Ga}}
\nc{\gpsmall}{\mbox{$\gp$}}
\nc{\gig}{\gimel}
\nc{\gns}{\sone(\Om,\gig)}
\nc{\nsr}[2]{#1}
\nc{\Srg}{{\mathbb{S}}}
\nc{\Srgs}{{\mathbb{S}^*}}
\nc{\NN}{{\w^{\w}}}
\nc{\ZN}{{\Z^{\bbN}}}
\nc{\NNup}{{\bbN^{\uparrow\bbN}}}
\nc{\Pof}{\op{P}}
\nc{\PN}{{\Pof(\w)}}
\nc{\roth}{{[\w]^{\w}}}
\nc{\Fin}{[\bbN]^{\text{$<\!\!\infty$}}} 
\nc{\ici}{[\bbN]^{ \infty, \infty}}
\nc{\Inc}{{\compactN^{\uparrow\bbN}}}
\nc{\powInc}[1]{{\big(\Inc\big)^{#1}}}
\nc{\powFin}[1]{{\big(\Fin\big)^{#1}}}
\nc{\powPN}[1]{{\big(\PN\big)^{#1}}}
\nc{\NcompactN}{{\compactN^\bbN}}
\nc{\Uarrow}{\smash{\big\uparrow}}
\nc{\LE}{\preccurlyeq}
\nc{\GE}{\succcurlyeq}
\nc{\op}{\operatorname}
\nc{\im}{\op{im}}
\nc{\Span}{\op{span}}
\nc{\maxfin}{\op{maxfin}}
\nc{\ran}{\op{range}}
\nc{\iso}{\cong}
\nc{\Madd}{{\M}^\star}
\nc{\cI}{\mathcal{I}}
\nc{\cJ}{\mathcal{J}}
\nc{\scrA}{\mathscr{A}}
\nc{\scrB}{\mathscr{B}}
\nc{\scrC}{\mathscr{C}}
\nc{\scrD}{\mathscr{D}}
\nc{\scrF}{\mathscr{F}}
\nc{\scrK}{\mathscr{K}}
\nc{\A}{\forall}
\nc{\B}{\mathrm{B}}
\nc{\cB}{\mathcal{B}}
\nc{\bB}{\mathbf{B}}
\nc{\BS}{\mathbf{B}(\mathcal{S})}
\nc{\BF}{\mathbf{B}(\mathcal{F})}
\nc{\BU}{\mathbf{B}(\mathcal{U})}
\nc{\cSp}{\mathcal{S}^+}
\nc{\cFp}{\mathcal{F}^+}
\nc{\cUp}{\mathcal{U}^+}

\nc{\BG}{\B_\Ga}
\nc{\BL}{\B_\Lambda}
\nc{\BT}{\B_\Tau}
\nc{\BTstar}{\B_{\Tau^*}}
\nc{\BO}{\B_\Om}
\nc{\DO}{\cD_\Om}
\nc{\KO}{\cK_\Om}
\nc{\CG}{C_\Ga}
\nc{\CL}{C_\Lambda}
\nc{\CT}{C_\Tau}
\nc{\CTstar}{C_{\Tau^*}}
\nc{\CO}{C_\Om}
\nc{\COgp}{C_{\Om^{\grbl}}}
\nc{\CLgp}{C_{\Lambda^{\grbl}}}
\nc{\BOgp}{\B_{\Om}^{\grbl}}
\nc{\BLgp}{\B_{\Lambda^{\grbl}}}
\nc{\sfC}{\mathsf{C}}
\nc{\sfD}{\mathsf{D}}
\nc{\bD}{\mathbf{D}}
\nc{\Tau}{\mathrm{T}}
\nc{\cA}{\mathcal{A}}
\nc{\cK}{\mathcal{K}}
\nc{\cD}{\mathcal{D}}
\nc{\cF}{\mathcal{F}}
\nc{\cS}{\mathcal{S}}
\nc{\cT}{\mathcal{T}}
\nc{\cG}{\mathcal{G}}
\nc{\cY}{\mathcal{Y}}
\nc{\J}{\mathcal{J}}
\nc{\cL}{\mathcal{L}}
\nc{\cM}{\mathcal{M}}
\nc{\cN}{\mathcal{N}}
\nc{\cH}{\mathcal{H}}
\nc{\cO}{\mathcal{O}}
\nc{\Op}{\mathrm{O}}
\nc{\rmA}{\mathrm{A}}
\nc{\rmF}{\mathrm{F}}
\nc{\rmB}{\mathrm{B}}
\nc{\rmD}{\mathrm{D}}
\nc{\rmP}{\mathrm{P}}
\nc{\cC}{\mathcal{C}}
\nc{\cP}{\mathcal{P}}
\nc{\bbQ}{\mathbb{Q}}
\nc{\bbR}{\mathbb{R}}
\nc{\cU}{\mathcal{U}}
\nc{\cQ}{\mathcal{Q}}
\nc{\Un}{\bigcup}
\nc{\cV}{\mathcal{V}}
\nc{\cW}{\mathcal{W}}
\nc{\Z}{{\mathbb Z}}
\nc{\Impl}{\Rightarrow}
\long\def\forget#1\forgotten{\marginpar{\textcolor{green}{Forgetting...}}}
\nc{\ft}{\mathfrak{t}}
\nc{\fb}{\mathfrak{b}}
\nc{\fc}{\mathfrak{c}}
\nc{\fd}{\mathfrak{d}}
\nc{\fg}{\mathfrak{g}}
\nc{\oo}{\infty}
\nc{\fr}{\mathfrak{r}}
\nc{\fk}{\mathfrak{k}}
\nc{\bidi}{\mathfrak{bidi}}
\nc{\fu}{\mathfrak{u}}
\nc{\fh}{\mathfrak{h}}
\nc{\fp}{\mathfrak{p}}
\nc{\fj}{\mathfrak{j}}
\nc{\fs}{\mathfrak{s}}
\nc{\w}{\omega}
\nc{\x}{\times}
\nc{\Iff}{\Leftrightarrow}

\nc{\nin}{\notin}
\nc{\cat}{\hat{\ }}
\nc{\sub}{\subseteq}
\nc{\spst}{\supseteq}
\nc{\sm}{\setminus}
\nc{\as}{\subseteq^*}
\nc{\les}{\le^*}
\nc{\leinf}{\le^{\infty}}
\nc{\leS}{\le_S}
\nc{\leF}{\le_{\mathcal{F}}}
\nc{\leU}{\le_{U}}
\nc{\rest}{\restriction}

\nc{\la}{\langle}
\nc{\ra}{\rangle}
\nc{\E}{\exists}
\nc{\dom}{\op{dom}}
\nc{\cov}{\op{cov}}
\nc{\add}{\op{add}}
\nc{\cof}{\op{cof}}
\nc{\cf}{\op{cf}}
\nc{\non}{\op{non}}
\nc{\unif}{\op{non}}
\nc{\COV}{\op{COV}}
\nc{\ADD}{\op{ADD}}
\nc{\COF}{\op{COF}}
\nc{\NON}{\op{NON}}
\nc{\impl}{\to}
\nc{\Lp}{\mathcal{L_\p}}
\nc{\Wlog}{without loss of generality}
\newtheorem{thm}{Theorem}[section]
\nc{\bthm}{\begin{thm}} \nc{\ethm}{\end{thm}}
\newtheorem{prop}[thm]{Proposition}
\nc{\bprp}{\begin{prop}} \nc{\eprp}{\end{prop}}
\newtheorem{fact}[thm]{Fact}
\nc{\bfct}{\begin{fact}} \nc{\efct}{\end{fact}}
\newtheorem{prob}[thm]{Problem}
\nc{\bprb}{\begin{prob}} \nc{\eprb}{\end{prob}}
\newtheorem{lem}[thm]{Lemma}
\nc{\blem}{\begin{lem}} \nc{\elem}{\end{lem}}
\newtheorem{claim}[thm]{Claim}
\nc{\bclm}{\begin{claim}} \nc{\eclm}{\end{claim}}
\newtheorem{cor}[thm]{Corollary}
\nc{\bcor}{\begin{cor}} \nc{\ecor}{\end{cor}}
\newtheorem{conj}[thm]{Conjecture}
\nc{\bcnj}{\begin{conj}} \nc{\ecnj}{\end{conj}}
\theoremstyle{definition}
\newtheorem{defn}[thm]{Definition}
\nc{\bdfn}{\begin{defn}} \nc{\edfn}{\end{defn}}
\newtheorem{obs}[thm]{Observation}
\nc{\bobs}{\begin{obs}} \nc{\eobs}{\end{obs}}
\theoremstyle{remark}
\newtheorem{rem}[thm]{Remark}
\nc{\brem}{\begin{rem}} \nc{\erem}{\end{rem}}
\newtheorem{cnv}[thm]{Convention}
\nc{\bcnv}{\begin{cnv}} \nc{\ecnv}{\end{cnv}}
\newtheorem{exam}[thm]{Example}
\nc{\bexm}{\begin{exam}} \nc{\eexm}{\end{exam}}
\nc{\bpf}{\begin{proof}} \nc{\epf}{\end{proof}}
\nc{\be}{\begin{enumerate}}
\nc{\ee}{\end{enumerate}}
\nc{\bi}{\begin{itemize}}
\nc{\bimy}{\my{\begin{itemize}}
\nc{\eimy}{\end{itemize}}}
\nc{\itm}{\item}
\nc{\ei}{\end{itemize}}
\nc{\Subsection}[1]{\goodbreak\subsection*{#1}}

\nc{\sone}{\mathsf{S}_1}
\nc{\sfin}{\mathsf{S}_\mathrm{fin}}
\nc{\ufin}{\mathsf{U}_\mathrm{fin}}
\nc{\Split}{\mathrm{Split}}

\nc{\gone}{\mathsf{G}_1}   
\nc{\Succ}{\mathrm{S}} 
 \nc{\gfin}{\mathsf{G}_\mathrm{fin}}

\DeclareMathOperator{\eexists}{\exists}
\DeclareMathOperator{\fforall}{\forall}
\nc{\Exists}[1]{\bigl(\eexists #1\bigr)}
\nc{\Forall}[1]{\fforall #1\ }
\nc{\Foralm}[1]{\fforall^* #1\ }
\nc{\End}[1]{#1}
\nc{\supp}{\op{supp}}
\nc{\bfP}{\mathbf{P}}

\DeclareMathOperator{\PONAG}{PONAG}
\nc{\Alice}{{\textsc{Alice}}}
\nc{\Bob}{{\textsc{Bob}}}
\DeclareMathOperator{\M}{M}

\title{Totally paracompact spaces and the Menger covering property}

\thanks{The research of the third author was funded by the National Science Center, Poland Weave-UNISONO call in the Weave programme
Project: Set-theoretic aspects of topological selections 2021/03/Y/ST1/00122 and also by the University of Messina
}

\author[M. Bonanzinga]{Maddalena Bonanzinga}
\address{Maddalena Bonanzinga, MIFT Department, University of Messina, Italy}
\email{mbonanzinga@unime.it}
\author[D. Giacopello]{Davide Giacopello}
\address{Davide Giacopello, MIFT Department, University of Messina, Italy}
\email{dagiacopello@unime.it}
\author[P. Szewczak]{Piotr Szewczak}
\address{Institute of Mathematics, Faculty of Mathematics and Natural Science,
College of Sciences, Cardinal Stefan Wyszy\'nski University in Warsaw, W\'oycickiego 1/3,
01–938 Warsaw, Poland}
\email{p.szewczak@wp.pl}
\urladdr{http://piotrszewczak.pl}

\subjclass[2010]{Primary: 54D20; 
Secondary: 03E17. 
}

\keywords{Total paracompactness, Menger property, GO-space, topological games}

\begin{document}

\maketitle

\begin{abstract}
A topological space is totally paracompact if any base of this space contains a locally finite subcover.
We focus on a problem of Curtis whether in the class of regular Lindel\"of spaces total paracompactness is equivalent to the Menger covering property.
To this end we consider topological spaces with certain dense subsets.
It follows from our results that the above equivalence holds in the class of Lindel\"of GO-spaces defined on subsets of reals.
We also provide a game-theoretical proof that any regular Menger space is totally paracompact and show that in the class of first-countable spaces the Menger game and a partial open neighborhood assignment game of Aurichi are equivalent.
We also show that if $\fb=\w_1$, then there is an uncountable subspace of the Sorgenfrey line whose all finite powers are Lindel\"of, which is a strengthening of a famous result due to Michael.
\end{abstract}

\section{Introduction}

By \emph{space} we mean a Hausdorff infinite topological space.
A space $X$ is \emph{totally paracompact}~\cite{Ford}, if every open base of $X$ contains a locally finite cover of $X$.
In 1924, Menger introduced the following covering property as an attempt to characterize $\sigma$-compactness in the class of metric spaces.
A metric space $X$ is \emph{Menger-basis} if for any basis $\cB$ of $X$, there are sets $B_0,B_1,\dotsc \in\cB$ whose diameters tend to zero and the family $\sset{B_n}{n\in\w}$ covers $X$.
After that Hurewicz reformulated the above definition, extending it to all topological spaces:
a space $X$ is \emph{Menger}, if for every sequence $\eseq{\cU}$ of open covers of $X$ there are finite sets $\cF_0\sub \cU_0, \cF_1\sub\cU_1,\dotsc$ such that the family $\bigcup\sset{\cF_n}{n\in \w}$ is a cover of $X$~\cite{H}. 
The Hurewicz definition of the Menger property is equivalent to the Menger-basis property in the class of metric spaces.
The Menger property is one of the central properties considered in the topological selections theory.

Our investigation is motivated by the following result of Curtis.

\bthm[{Curtis \cite[Theoirem~1.1.]{Curtis}}]\label{thm:Curtis} Every regular Menger space is totally paracompact. 
\ethm

Notice that there are non-regular countable spaces that may not be paracompact, e.g., Bing's space~\cite{SteenSeebach}.

A space is \emph{Lindel\"of} if any open cover of this space contains a countable subcover (notice that in this approach, we do not request that the space is regular).
Since any Menger space is Lindel\"of, in the light of the above result, Curtis formulated the following problem.

\bprb[{Curtis~\cite{Curtis}}]\label{Q:main}
Is it true that any Lindel\"of totally paracompact space is Menger?
\eprb

Thus far it is known that there is a positive solution to Problem~\ref{Q:main} in the class of metrizable spaces~\cite{lelek} and in the class of generalized ordered spaces defined on subsets of the real line~\cite{Sakai}.
In Section~\ref{sec:totparmen}, we show that total paracompactness is equivalent to the Menger property in the class of Lindel\"of spaces having $\sigma$-compact dense subset.
We also consider another classes of spaces, where it holds.
Our results capture all known classes of spaces, where the Curtis conjecture is true.
Since the Menger property can be characterized using a topological game~\cite{H, SzTsGame}, in Section~\ref{sec:curtisgame}, we provide a game-theoretical proof of Theorem~\ref{thm:Curtis} which seems to be much simpler than the original one (modulo game-characterization of the Menger property).
Section~\ref{sec:aurichi} is devoted to relations between a partial open neighborhood assignment game of Aurichi~\cite{Aurichi} and the Menger game (see definitions below).
In particular, we show that in the class of first-countable spaces these two games are equivalent.
The Aurichi game comes on stage naturally when considering Menger spaces and totally paracompact spaces in the context of so-called \emph{$D$-spaces} and an old problem whether every regular Lindel\"of space is a $D$-space~\cite{D} (see Section~5 for definitions and more details).
We finish our paper with comments and open problems.
\section{Menger vs. total paracompactness}\label{sec:totparmen} 

We establish the following partial answer to Problem~\ref{Q:main}.
Recall that a paracompact space with a dense $\sigma$-compact subspace is Lindel\"of~\cite[Theorem~5.1.25]{engelking}.

\bthm\label{Th:tot_par_impl_Men1}
Every totally paracompact space with a $\sigma$-compact dense subset is Menger.
\ethm

\bpf
Assume that $X$ is a totally paracompact space with a $\sigma$-compact dense subset $\bigcup_{n\in \omega} C_n$, where each set $C_n$ is compact.
Fix a sequence $\eseq{\cU}$ of open covers of $X$.
We may assume that the families $\cU_n$ are closed under finite unions.
Fix $x\in X$.
For each $n\in \omega$, let
$\cB_n(x)$ be the family of all open neighborhoods $V$ of $x$ such that $V$ is contained in some set from $\cU_n$ and $V\cap C_n\neq\emptyset$.
Then the collection $\cB(x):=\bigcup_{n \in \omega}\cB_n(x)$ is a local base at $x$.
It follows that the family $\cB:=\bigcup_{x\in X}\cB(x)$ is a base for $X$.
By the assumption, there exists a locally finite cover ${\cB}^{\prime} \sub {\cB}$ of $X$.
Fix $n \in \omega$.
The family $\cB'$ is locally finite and it covers the compact set $C_n$.
Thus, the set
\[
\sset{B\in\cB'}{B\cap C_n\neq\emptyset}
\]
is finite and the set
\[
\cB_n':=\cB'\cap \bigcup_{x\in X}\cB_n(x)
\]
is finite, too.
By the definition, each element from the union $\bigcup_{x\in X}{\cB}_n(x)$ is contained in a set from $\cU_n$.
Since the family $\cU_n$ is closed under finite unions, there is a set $U_n\in\cU_n$ such that $\bigcup{\cB}_n'\sub U_n$.
Since the family $\cB'=\Un_{n\in\w}\cB'_n$ covers $X$, the family $\{U_n: n\in\omega\}$ covers $X$, too.
\epf

\bcor\label{Cor:tot_par_impl_Men1}
Every separable totally paracompact space is Menger.
\ecor

\bcor[Lelek~\cite{lelek}]
Every metrizable totally paracompact space is Menger.
\ecor

A space $X$ is \emph{totally metacompact}~\cite{telg}, if every open base of $X$ contains a point-finite cover of $X$. Clearly, every totally paracompact space is totally metacompact. We asked if we could replace total paracompactness in the statement of Theorem~\ref{Th:tot_par_impl_Men1} with total metacompactness.
The following example gives a negative answer to this question.

\bexm
There is a totally metacompact, not paracompact, not Lindel\"of space having a $\sigma$-compact dense subset:
Consider $X:=(\w+1)\x(\w_1+1)\sm \{(\w,\w_1)\}$ with a topology inherited from the Tychonoff product of the one-point compactification of $\w$ and the one-point compactification of the discrete space $\w_1$.
Such a space is not normal, and thus it is neither paracompact nor Lindel\"of.
The set $\w\times (\w_1+1)$ is a $\sigma$-compact dense subset of $X$.
The space $X$ is totally metacompact~\cite{B}.
\eexm

Recall that every metacompact separable space is Lindel\"of. Then we found another partial answer to the Curtis question with the following result.
The proof of it is a modification of the proof of Theorem~\ref{Th:tot_par_impl_Men1}.
 
\bthm \label{thm:sep}
Every separable totally metacompact space is Menger.
\ethm
\bpf
Let $X$ be a totally metacompact space with a countable dense subset $Q=\sset{q_n}{n\in\w}$.
Fix a sequence $\eseq{\cU}$ of open covers of $X$.
We may assume that each family $\cU_n$ is closed under finite unions.
Fix $x \in X$.
For each $n \in \omega$, let $\cB_n(x)$ be the family of all open sets such that $\{x, q_n\} \sub V$ and $V$ is contained in some set from $\cU_n$.
Then the collection $\cB(x):=\bigcup_{n \in \w}\cB_n(x)$ is a local base at $x$.
It follows that the family $\cB:=\bigcup_{x\in X}\cB(x)$ is a base for $X$.
By the assumption, there exists a point-finite cover $\cB'\sub\cB$ of $X$.
Fix $n\in\w$.
Since the family $\cB'$ is point-finite, the set
\[
\sset{B\in\cB'}{q_n\in B}
\]
is finite and the set
\[
\cB_n':=\cB'\cap \Un_{x\in X}\cB_n(x)
\]
is finite, too.
By the definition, each set from the union $\Un_{x\in X}\cB_n(x)$ is contained in a set from $\cU_n$.
Since the family $\cU_n$ is closed under finite unions, there is a set $U_n\in\cU_n$ such that $\bigcup{\cB}_n'\sub U_n$.
Since the family $\cB'=\Un_{n\in\w}\cB'_n$ covers $X$, the family $\sset{U_n}{n\in\w}$ covers $X$, too.
\epf

\bcor
For a separable space $X$, the following assertions are equivalent.
\be
\item The space $X$ is regular and Menger.
\item The space $X$ is totally paracompact.
\item The space $X$ is regular and totally metacompact. 
\ee
\ecor

A space $X$ is \emph{scattered} if for every nonempty closed set $A\sub X$ there is a point $a\in A$ which is isolated in the relative topology of $A$.
Let $X$ be a space.
Define $X^{(0)}:=X$.
Fix an ordinal number $\alpha$ and assume that the sets $X^{(\beta)}$ have been defined for all ordinal numbers $\beta<\alpha$.
If $\alpha=\beta+1$, for some $\beta$, then
\[
X^{(\alpha)}:=X^{(\beta)}\sm\sset{x\in X^{(\beta)}}{x\text{ is isolated in }X^{(\beta)}}.
\]
If $\alpha$ is a limit ordinal, then
\[
X^{(\alpha)}=\bigcap_{\beta<\alpha}X^{(\beta)}.
\]
A space $X$ is scattered if and only if, there is an ordinal number $\alpha$ such that $X^{(\alpha)}=\emptyset$.

\bprp\label{prp:scatt} 
Let $X$ be a Lindel\"of space such that the space $X^{(\alpha)}$ is Menger for some ordinal number $\alpha$.
Then the space $X$ is Menger.
\eprp

\bpf
Let $X$ be a Lindel\"of scattered space and $\alpha$ be the minimal ordinal number such that $X^{(\alpha)}$ is Menger.
The proof is by induction on $\alpha$.
If $\alpha=0$, then there is nothing to prove.
Fix an ordinal number $\alpha$ and assume that for each Lindel\"of space $Y$, if $Y^{(\beta)}$ is Menger, for some $\beta<\alpha$, then the space $Y$ is Menger.
Fix a sequence $\eseq{\cU}$ of open covers of $X$.
Let
\[
X':=X\sm X^{(\alpha)}.
\]
There are finite sets $\cF_0\sub\cU_0, \cF_1\sub\cU_1,\dotsc$ such that the family $\Un_{n\in\w}\cF_n$ covers $X^{(\alpha)}$.
Let $U:=\Un\Un_{n\in\w}\cF_n$.
Fix $x\in X\sm U$.
Then there is an open neighborhood $V_x$ of $x$ in $X$ such that $\overline{V_x}^{(\beta)}=\emptyset$.
Thus, the set $\overline{V_x}$ is Menger, by the assumption.
Since the space $X\sm U$ is Lindel\"of as a closed subspace of $X$, there is a countable set $C\sub X\sm U$ such that 
\[
X\sm U\sub \Un_{x\in C}\overline{V_x}.
\]
Then the space $X\sm U$ is Menger as a countable union of Menger spaces $\overline{V_x}\sm U$ for $x\in C$.
There are finite sets $\cF'_0\sub\cU_0,\cF'_1\sub\cU_1,\dotsc$ such that the family $\Un_{n\in\w}\cF'_n$ covers $X\sm U$.
Consequently the family $\Un_{n\in\w}\cF_n\cup\cF_n'$ covers $X$ and the sets $\cF_n\cup\cF_n'$ are finite subsets of $\cU_n$.
\epf

\bcor\label{scattered}
Every totally metacompact regular Lindel\"of space $X$ such that  $X^{(\alpha)}$ is separable, for some ordinal number $\alpha$, is Menger.
\ecor

\bpf
If $X$ is totally metacompact, then for every ordinal $\alpha$, the set $X^{(\alpha)}$ is also totally metacompact.
By Theorem~\ref{thm:sep}, the space $X^{(\alpha)}$ is Menger.
Apply Proposition~\ref{prp:scatt}.
\epf

\section{GO-spaces}\label{sec:GO}

Sakai~\cite{Sakai} considered Problem~\ref{Q:main} for subspaces of the Sorgenfrey line and more generally in the class of generalized ordered spaces which are defined on subsets of the real line (see definition below).
This class of spaces has a long history in general topology and it has been well studied.
A \emph{real generalized-ordered space} (or abbreviated \emph{real GO-space}) is a triple $(X,<,\tau)$, where $X\sub\bbR$, the relation $<$ is the usual order on $\bbR$ restricted to $X$ and $\tau$ is a Hausdorff topology on $X$ generated by some intervals (possibly degenerated).
For simplicity we write a real GO-space $(X,<,\tau)$ as $X$.
Every real GO-space $X$ we can split 
\[
X= R_X\cup E_X\cup L_X\cup D_X,
\]
where:
\begin{align*}
D_X & :=  \sset{x\in X}{\{x\}\text{ is open in }X},\\
R_X  & :=  \sset{x\in X\sm D_X}{[x,\infty) \text{ is open in  }X},\\
L_X  & := \sset{x\in X\sm D_X}{(-\infty,x]\text{ is open in  }X},\\
E_X  & :=  \sset{x\in X\sm D_X}{\sset{(a,b)}{a,b\in X, x\in (a,b)}\text{ is a base at }x}.
\end{align*}

We get the following corollary from the previous section.

\bcor\label{GO}
Every totally metacompact Lindel\"of real GO-space is Menger.
\ecor

\bpf There is an ordinal number $\alpha$ such that $X^{(\alpha)}=X^{(\alpha+1)}$.
It follows that $X^{(\alpha)}$ is a closed subspace of $X$ with no isolated points.
Any real GO-space with no isolated points is separable.
Apply Theorem~\ref{scattered}.
\epf

\bcor\label{cor:sakai}
For a Lindel\"of real GO-space $X$, the following assertions are equivalent.
\be
\item The space $X$ is Menger.
\item The space $X$ is totally paracompact.
\item The space $X$ is totally metacompact.
\ee
\ecor

\brem
The above corollary was obtained also by Sakai~\cite[Theorem 5.2]{Sakai}, but his proof is much more technical than our, more general approach.
\erem

Now we consider real GO-spaces defined on the entire real line. So, let $X=\Bbb R$, $X= R_X\cup E_X\cup L_X\cup I_X$, and let $\tau$ be some topology 
$\tau$ that comes from the partition. For a subset $C\subseteq X$,  define
\begin{align*}
i(C) & :=  \sset{x\in X}{\{x\}\text{ is isolated in } \tau|C},\\
r(C) & :=  (R_X\cap C)\sm i(C),\\
l(C) & :=  (L_X\cap C)\sm i(C),\\
e(C) & :=  E_X\cap C.
\end{align*}

Recall the following result of Benett--Balogh.

\bthm[{Benett--Balogh\cite{BB}}]\label{BB}
Let $X$ be a real GO-space defined on the entire real line. The space $X$ is totally paracompact if and only if for every perfect set $C$ in the standard topology on the real line, the set $i(C)$ is not dense and co-dense in $C$ in the standard topology.
\ethm

\bthm 
A real GO-space defined on the entire real line is Menger if and only if it is $\sigma$-compact.
\ethm
\bpf
Every $\sigma$-compact space is Menger.
So, assume that $X$ is a Menger real GO-space defined on the entire real line. 
Let ${\cU}$ be the collection of all intervals $(a,b)$, where $a,b\in\bbR$ and $a<b$, such that $(a,b)\cap X$ is $\sigma$-compact in $X$.
Then the set $U:=\bigcup \cU$ is $\sigma$-compact in $X$.
Towards a contradiction, assume that the set $C:=X\sm U$ is nonempty. 
Since $U$ is open in the standard topology, we have that the set $C$ is closed in the standard topology.
Now we prove that $C$ is perfect.
Assume that there is an interval $(a,b)$ such that $(a,b)\cap C=\{x\}$ for some $x\in\bbR$.
Since $(a,x)\sub U$, the set $U$ is $\sigma$-compact in $X$ and $(a,x)$ is an $F_\sigma$-set, the set $(a,x)$ is $\sigma$-compact in $X$.
Similarly, the set $(x,b)$ is $\sigma$-compact in $X$.
Thus, the set $(a,b)\cap X$ is $\sigma$-compact in $X$, too.
It follows that $(a,b)\sub U$, a contradiction.

The set $C\cap R_X$ is not dense in $C$ with the standard topology:
Suppose not.
Since $C$ is perfect, there is a copy $C^\prime$ of the Cantor set inside $C$ such that the sets $i(
C')\sub R_X\cap C$ and $r(C')\sub R_X\cap C$ are dense in $C'$ with respect to the standard topology on $C'$.
Then the set $i(C')$ of isolated points of $C^\prime$ is dense and co-dense in $C^\prime$.
By Theorem~\ref{BB} and~Corollary~\ref{cor:sakai} we get a contradiction.
Similarly, the set $C\cap L_X$ is not dense in $C$ with the standard topology.

By the above observations, there is an interval $[a,b]$ in $\bbR$ such that the set $C':=[a,b]\cap C$ is perfect in the standard topology and it is disjoint with $R_X\cup L_X$.
Thus, we may assume that $C=C'$.
Now we prove that the set $i(C)$ is dense in $C$ with respect to the standard topology:
Suppose not and let $[a,b]$ be an interval in $\bbR$ such that the set $[a,b]\cap C$ is perfect in the standard topology and it is disjoint with $i(C)$.
Then the topologies on $[a,b]\cap C$ inherited from $\bbR$ and $X$ are the same.
It follows that $(a,b)\cap X$ is $\sigma$-compact, and thus $(a,b)\in\cU$, a contradiction.

The set $i(C)$ is also co-dense in $C$ with respect to the standard topology:
If not, then there is an interval $[a,b]$ in $\bbR$ such that the set $C':=[a,b]\cap C$ is perfect in the standard topology and it is contained in $i(C)$.
It follows that $C'$ is a closed and discrete subspace of $X$.
Since $X$ is Lindel\"of, the set $C'$ is countable, a contradiction.

Finally, we have that the set $i(C)$ is dense and co-dense in $C$ in the standard topology.
By Theorem~\ref{BB} and~Corollary~\ref{cor:sakai} we get a contradiction.

We conclude that $C=\emptyset$, and thus the space $X$ is $\sigma$-compact.
\epf


\section{Curtis' result via topological games}
\label{sec:curtisgame}

Many classical covering properties, starting with the Menger property, admit natural game-theoretic counterpart \cite{Telg1, SzTsGame, ABG}.
Recall the following definition.

\bdfn
A \emph{Menger game} played on a space $X$, denoted by $\M(X)$, is a game with two players, \Alice{} and \Bob{}, with an inning per each natural number $n$.
In round $n$ \Alice{} plays an open cover ${\cU}_n$ of $X$ and \Bob{} replies with a finite set $\cF_n\sub \cU_n$.
\Bob{} wins the game, if the family $\Un_{n\in\w}\cF_n$ covers $X$ and \Alice{} wins otherwise.
\edfn

\bthm[{\cite{H, SzTsGame}}]\label{thm:M}
A space $X$ is Menger if and only if \Alice{} has no winning strategy in the game $\M(X)$.
\ethm

We provide a game-theoretical proof of Theorem~\ref{thm:Curtis} which seems to be much simpler than the original one (modulo the above game-characterization of the Menger property).
We need the following straightforward observation.

\bobs\label{obs:str}
Let ${\cB}^\prime$ be an open cover of a space $X$ and $A^\prime$ be a closed subset of $X$.
If 
\bi
\item $A^\prime\sub B^\prime$ for some $B^\prime\in{\cB}^{\prime}$,
\item $B\cap A^\prime=\emptyset$ for every $B\in{\cB}^\prime\setminus\{B^\prime\}$, and
\item $\sset{V_B}{B\in{\cB}^\prime}$ is an open cover of $X$,
\ei
then $A^\prime\sub V_{B^\prime}$.
\eobs

\bthm[{\cite[Theorem~1.1.]{Curtis}}]
Every regular Menger space is totally paracompact.
\ethm
\bpf
Let $X$ be a regular Menger space and $\cB$ be a fixed base for $X$.

\textbf{Round 0:} By Lindel\"ofness, there exists a countable subcover, say ${\cB}_0$, of $\cB$. Again, by Lindel\"ofness, there is an open cover $\sset{V^0_B}{B\in \cB_0}$ such that $\overline{V^0_B}\sub B$ for all $B\in \cB_0$.
\Alice{} plays $\cV_0:=\sset{V^0_B}{B\in \cB_0}$.
\Bob{} replies with s $\{V^0_B: B\in{\cF}_0\}$, where ${\cF}_0$ is a finite subset of ${\cB}_0$.
Put $A_0:=\bigcup_{B\in{\cF}_0}\overline{V_B^0}$.

\textbf{Round 1:}
There is a countable subcover ${\cB}_1$ of ${\cB}$ such that ${\cF}_0\subseteq{\cB}_1$ and for every set $B\in{\cB}_1\sm{\cF}_0$, we have $B\cap A_0=\emptyset$.
There is an open cover $\{V_B^1: B\in{\cB}_1\}$ such that $\overline{V_B^1}\sub B$ for all $B\in\cB_1$.
Put $V_1:=\bigcup_{B\in{\cF}_0}V_B^1$.
By Observation~\ref{obs:str}, we have that $A_0\subseteq V_1 \subseteq \overline{V_1}\subseteq \bigcup{\cF}_0$. 
\Alice{} plays ${\cV}_1:=\{V_1\cup V_B^1: B\in{\cB}_1\setminus{\cF}_0\}$.
\Bob{} replies with $\sset{V_1\cup V_B^1}{B\in{\cF}_1\sm {\cF}_0}\sub{\cV}_1$, where ${\cF}_1$ is a finite subset of ${\cB}_1$.
We may assume that ${\cF}_0\subseteq{\cF}_1$.
Put 
\[
A_1:=\bigcup_{B\in{\cF}_1\setminus {\cF}_0}\overline{V_1\cup V_B^1}=\bigcup_{B\in{\cF}_1}\overline{V_B^1}.
\]
Fix a natural number $n>1$.

\textbf{Round $\mathbf{n}$:}
There is a countable subcover ${\cB}_n$ of ${\cB}$ such that ${\cF}_{n-1}\subseteq{\cB}_n$ and for every set $B\in{\cB}_n\sm{\cF}_{n-1}$, we have $B\cap A_{n-1}=\emptyset$.
There is an open cover $\{V_B^n: B\in{\cB}_n\}$ such that $\overline{V_B^n}\sub B$ for all $B\in\cB_n$.
Put $V_n:=\bigcup_{B\in{\cF}_{n-1}}V_B^n$.
By Observation~\ref{obs:str}, we have that $A_{n-1}\subseteq V_n \subseteq \overline{V_n}\subseteq \bigcup{\cF}_{n-1}$. 
\Alice{} plays ${\cV}_n:=\{V_n\cup V_B^n: B\in{\cB}_n\setminus{\cF}_{n-1}\}$.
\Bob{} replies with a family $\sset{V_n\cup V_B^n}{B\in{\cF}_n\sm {\cF}_{n-1}}\sub{\cV}_n$, where ${\cF}_n$ is a finite subset of ${\cB}_n$.
We may assume that ${\cF}_{n-1}\subseteq{\cF}_n$.
Put 
\[
A_n:=\bigcup_{B\in{\cF}_n\setminus {\cF}_{n-1}}\overline{V_n\cup V_B^n}=\bigcup_{B\in{\cF}_n}\overline{V_B^n}.
\]

Since \Alice{} has no winning strategy in $\M(X)$, there is a play, where Alice uses the above defined strategy and the game is won by \Bob{}, i.e., 
\[
\bigcup\{V_B^n: B\in{\cF}_n, n\in\omega\}
\]
is an open cover of $X$.
For every $n\in\w$ and $B\in\cF_n$ we have $V^n_B\sub B$.
It follows that $\bigcup_{n\in\omega}{\cF}_n\subseteq{\cB}$ is an open cover of $X$.
Fix $x\in X$.
There is a natural number $n$ such that $x\in V_n$.
For $k>n$ and $B\in \cF_k\sm \cF_n$ we have $B\cap V_n=\emptyset$.
It follows that the family $\Un_{n\in\w}\cF_n$ is locally finite.
Thus, the space $X$ is totally paracompact.
\epf

\section{Partial open neighborhood assignment game of Aurichi}
\label{sec:aurichi}

Let $X$ be a space.
An \emph{open neighborhood assignment of $X$} is a family $\sset{U_x}{x\in X}$ of open sets in $X$ such that $x\in U_x$ for all $x\in X$.
The space $X$ is a \emph{$D$-space}~\cite{D} if for any open neighborhood assignment $\sset{U_x}{x\in X}$ of $X$, there is a closed and discrete in $X$ space $D\sub X$ such that the family $\sset{U_x}{x\in D}$ covers $X$.
In the literature, many references to a result that any Menger space is a $D$-space lead to a paper of Aurichi~\cite{Aurichi}.
It turns out that this result, in the class of regular spaces, follows directly from much older Theorem~\ref{thm:Curtis} (Aurichi does not assume that a space is regular).
We provide a proof for the sake of completeness.

\bprp\label{prop:tot_par_impl_D}
	Every totally paracompact space is a D-space.
\eprp
\bpf
Let $X$ be a space and $\{U_x:x\in X\}$ be an open neighborhood assignment of $X$.
For $x\in X$, let ${\cB}_x$ be a local base at $x$ such that $\bigcup{\cB}_x\subseteq U_x$. 
Then ${\cB}=\bigcup_{x\in X}{\cB}_x$ is a base for $X$.
By the assumption, there is a family ${\cB}^\prime \subset {\cB}$ which is a locally finite cover of $X$.
Clearly, the sets ${\cB}^\prime \cap {\cB}_x$ are finite for each $x\in X$.
Then set $D=\{x\in X: {\cB}^\prime \cap {\cB}_x\not=\emptyset\}$ is closed and discrete in $X$ and the family $\{U_x: x\in D\}$ covers $X$.
\epf

In order to prove that any Menger space is a $D$-space, Aurichi used the following topological game.
Let $X$ be a space.
A \emph{partial open neighborhood assignment in} $X$ is a family $\sset{V_x}{x\in Y}$, where $Y\sub X$ and each set $V_x$ is an open neighborhood of $x$.
For a partial open neighborhood assignment $\sset{V_x}{x\in Y}$ in $X$ and a set $D\sub Y$ define $V_D:=\Un_{x\in D}V_x$.  
Partial open neighborhood assignments $\sset{V_x}{x\in Y}$ and $\sset{V'_x}{x\in Y'}$ in $X$ are \emph{compatible} if $V_x=V'_x$ for each $x\in Y\cap Y'$.
Since the ambient space $X$ will be clear from the context, a partial open neighborhood assignment in $X$, we just call a \emph{partial open neighborhood assignment}

\bdfn[{Aurichi~\cite[Definition~2.1]{Aurichi}}]
Let $X$ be a space.
A \emph{partial open neighborhood assignment game} played on $X$, denoted by $\PONAG(X)$, is a game with two players, \Alice{} and \Bob{}, with an inning per each natural number $n$.
In round $0$, \Alice{} plays a partial open neighborhood assignment $\cU_0=\sset{V_x}{x\in Y_0}$ for some $Y_0\subseteq X$ such that $\cU_0$ covers $X$.
Then \Bob{} replies with a closed and discrete in $X$ space $D_0\sub Y_0$.
In round $n>0$, \Alice{} plays a partial open neighborhood assignment $\cU_n=\sset{V_x}{x\in Y_n}$ in $X$ such that
\bi
\item $Y_n\sub X\sm \Un_{k<n}V_{D_k}$,
\item $\cU_n$ covers $X\sm \Un_{k<n}V_{D_k}$.
\ei
\Bob{} replies with a closed and discrete in $X$ space $D_n\sub Y_n$.

\Bob{} wins the game, if the family $\Un_{n\in\w}\sset{V_x}{x\in D_n}$ covers $X$, otherwise \Alice{} wins.
\edfn

Recall some important results from the Aurichi paper from our point of view.

\bthm\label{thm:M-PONAG}\cite{Aurichi}
If $X$ is a Menger space, then \Alice{} has no winning strategy in the game $\PONAG(X)$.
\ethm

\bthm \cite{Aurichi}
Let $X$ be a space.
If \Alice{} has no winning strategy in the game $\PONAG(X)$, then $X$ is a $D$-space.
\ethm

It turns out that in a wide class of topological spaces, the implication from Theorem~\ref{thm:M-PONAG} is reversible.

\bthm \label{Alice_game_equivalence}
Let $X$ be a Lindel\"of space such that each closed non-Menger subspace of $X$ contains a convergent sequence.
\Alice{} has a winning strategy in the game $\PONAG(X)$ if and only if \Alice{} has a winning strategy in the game $\M(X)$.
\ethm

\bpf
($\Rightarrow$): 
By Theorem~\ref{thm:M-PONAG}, the space $X$ is not Menger.
Apply Theorem~\ref{thm:M}.

($\Leftarrow$):
Fix a winning strategy $\sigma$ for \Alice{} in the game $\M(X)$. 
Play the game $\M(X)$ according to $\sigma$ in order to show that we may assume some additional properties of moves for both players.

\textbf{Round 0:} Let $\cU_0$ be a move of \Alice{} played according to the strategy $\sigma$.
Since the space $X$ is Lindel\"of, we may assume that $\cU_0=\sset{U^0_k}{k\in\w}$ is an increasing cover of $X$.
By the assumption, there is a convergent sequence $Y_0=\sset{y^0_k}{k\in\w}\sub X$. 
Since the cover $\cU_0$ is increasing, a cofinal subset of $Y_0$ is contained in a single set from $\cU_0$.
Thus, we may assume that $Y_0\sub U^0_0$.
If not, then it is enough to replace $Y_0$ by its cofinal subset and $\cU_0$ by its cofinal subfamily.
Then, we may assume that \Bob{} replies with a single set $U_0\in\cU_0$.

\textbf{Round 1:}
Let $\cU_1$ be the next move of \Alice{} according to $\sigma$.
We may assume that $\cU_1=\sset{U^1_k}{k\in\w}$ is an increasing open cover of $X$ such that $U_0\sub U^1_0$.
Since \Alice{} has a winning strategy in $\M(X)$, the space $X\sm U_0$ cannot be Menger.
Indeed, the strategy $\sigma$ indicates a strategy for \Alice{} in $\M(X\sm U_0)$.
If $X\sm U_0$ is Menger, then by Theorem~\ref{thm:M} this indicated strategy cannot be a winning one which causes that $\sigma$ is not a winning strategy, a contradiction.
By the assumption, the space $X\sm U_0$ contains a convergent sequence $Y_1=\sset{y^1_k}{k\in\w}$.
Since the cover $\cU_1$ is increasing, a cofinal subset of $Y_1$ is contained in a single set from $\cU_1$.
Thus, we may assume that $Y_1\sub U^1_0\sm U_0$.
If not, then it is enough to replace $Y_1$ by its cofinal subset and $\cU_1$ by its cofinal subfamily.

Fix $n>0$ and let $U_{n-1}$ be the last set picked by \Bob{} in round $n-1$ in $\M(X)$.

\textbf{Round $\mathbf{n}$: }
Let $\cU_n$ be a family played by \Alice{}, according to $\sigma$.
In a similar manner as in round 1, we may assume that $\cU_n=\sset{U^n_k}{k\in\w}$ is an increasing open cover of $X$ such that $U_{n-1}\sub U^n_0$ and there is a convergent sequence $Y_n\sub U^n_0\sm U_{n-1}$.
Then \Bob{} replies with a set $U_n\in\cU_n$.

Formally, if \Alice{} has a winning strategy in the game $\M(X)$, then \Alice{} has a winning strategy in $\M(X)$ satisfying all the above additional conditions.
Thus, we assume that the strategy $\sigma$ already has all these properties.
If 
\[
(\cU_0,U_0,\cU_1,U_1,\dotsc)
\]
is a play in the game $\M(X)$, where \Alice{} uses the strategy $\sigma$ and in each inning \Bob{} picks a single set, as above, then the family $\sset{U_n}{n\in\w}$ is increasing and it does not cover $X$.

Now we translate the above strategy $\sigma$ to a strategy for \Alice{} in $\PONAG(X)$ as follows.

\textbf{Round 0:}
Let $\cU_0=\sset{U^0_k}{k\in\w}$ be the family played by \Alice{} according to $\sigma$ in $\M(X)$ and the set $Y_0=\sset{y^0_k}{k\in\w}\sub U^0_0$ be a convergent sequence as above.
Let $V_{y^0_k}:=U^0_k$ for $k\in\w$.
In $\PONAG(X)$, \Alice{} plays a partial open neighborhood assignment $\sset{V_{y^0_k}}{y^0_k\in Y_0}$ and \Bob{} replies with a closed and discrete subspace $D_0\sub Y_0$ in $X$.
Since $Y_0$ is a convergent sequence, the set $D_0$ is finite.
Since the family $\cU_0$ is increasing, there is a set $U_0\in\cU_0$ such that $V_{D_0}=U_0$.
In $\M(X)$, \Bob{} replies with $U_0$.

\textbf{Round 1:}
Let $\cU_1$ be the next move of \Alice{} according to $\sigma$ in $\M(X)$ and $Y_1\sub U^1_0\sm U_0$ be convergent sequence as above.
Let $V_{y^1_k}:=U^1_k$ for $k\in\w$.
In $\PONAG(X)$, \Alice{} plays a partial open neighborhood assignment $\sset{V_{y^1_k}}{y^1_k\in Y_1}$ and \Bob{} replies with a closed and discrete subspace $D_1\sub Y_1$ in $X$.
Since $Y_1$ is a convergent sequence, the set $D_1$ is finite.
Since the family $\cU_1$ is increasing, there is a set $U_1\in\cU_1$ such that $V_{D_1}=U_1$.
In $\M(X)$, \Bob{} replies with $U_1$.

Fix $n>0$ and suppose that open sets $U_0\sub \dotsb\sub U_{n-1}$, partial open neighborhood assigments $\sset{V_x}{x\in Y_k}$ and finite sets $D_k\sub Y_k\sub U_k\sm U_{k-1}$ have already been defined for $0\leq k\leq n$.

\textbf{Round $\mathbf{n}$:} Let $\cU_n$ be a family played by \Alice{}, according to $\sigma$ in $\M(X)$ and $Y_n=\sset{y^n_k}{k\in\w}\sub U^n_0\sm U_{n-1}$ be a convergent sequence as above.
Let $V_{y^n_k}:=U^n_k$ for $k\in\w$.
In $\PONAG(X)$, \Alice{} plays a partial open neighborhood assignment $\sset{V_{y^n_k}}{y^n_k\in Y_n}$ and \Bob{} replies with a closed and discrete subspace $D_n\sub Y_n$ in $X$.
Since $Y_n$ is a convergent sequence, the set $D_n$ is finite.
Since the family $\cU_n$ is increasing, there is a set $U_n\in\cU_n$ such that $V_{D_n}=U_n$.
In $\M(X)$, \Bob{} replies with $U_n$.

Since $\sigma$ is a winning strategy for \Alice{} in $\M(X)$, the family $\sset{U_n}{n\in\w}$ does not cover~$X$.
We have 
\[
\Un_{n\in\w}U_n=\Un_{n\in\w}V_{D_n},
\]
and thus the above defined strategy for \Alice{} in $\PONAG(X)$ is a winning strategy.
\epf

\bprp
Let $X$ be a space.
If \Bob{} has a winning strategy in $\M(X)$, then \Bob{} has a winning strategy in $\PONAG(X)$. 
\eprp

\bpf
Let $\sigma$ be a winning strategy for \Bob{} in $\M(X)$.
Define a strategy for \Bob{} in $\PONAG(X)$ as follows.

\textbf{Round 0:}
Assume that in $\PONAG(X)$, \Alice{} plays a partial open neighborhood assignment $\cU_0=\sset{V^0_x}{x\in Y_0}$ for some set $Y_0\sub X$.
A family $\cU_0':=\cU_0$ is an open cover of $X$.
Assume that $\cU_0'$ is a family played by \Alice{} in $\M(X)$.
Let $\cF_0'\sub\cU'_0$ be a finite set played by \Bob{} in $\M(X)$, according to $\sigma$, as a response to $\cU_0'$.
There is a finite set $D_0\sub Y_0$ such that $\Un\cF'_0=V_{D_0}$.
Then in $\PONAG(X)$, \Bob{} plays the set $D_0$.

\textbf{Round 1:}
Assume that, in $\PONAG(X)$, \Alice{} plays an open neighborhood assignment $\cU_1=\sset{V^1_x}{x\in Y_1}$ for some set $Y_1\sub X\sm V_{D_0}$.
Since the family $\cU_1$ covers $X\sm V_{D_0}$, a family $\cU_1':=\sset{V_x^1\cup V_{D_0}}{x\in Y_1}$ is an open cover of $X$.
Assume that $\cU_1'$ is the next move of \Alice{} in $\M(X)$.
Let $\cF_1'\sub\cU_1'$ be a finite set played by \Bob{} in $\M(X)$, according to $\sigma$, as a response to $\cU_1'$.
There is a finite set $D_1\sub Y_1$ such that $\cF'_1=\sset{V^1_x\cup V_{D_0}}{x\in D_1}$.
We have $\Un\cF_1'=V_{D_0}\cup V_{D_1}$.
Then in $\PONAG(X)$, \Bob{} plays the set $D_1$.

Fix $n>1$ and assume that open neighborhoods assignments $\sset{V_x}{x\in Y_k}$ played by \Alice{} and sets $D_k\sub Y_k$ played by \Bob{} have been already defined for all $0\leq k\leq n$.

\textbf{Round $\textbf{n}$:}
Assume that, in $\PONAG(X)$, \Alice{} plays an open neighborhood assignment $\cU_n=\sset{V^n_x}{x\in Y_n}$ in $\PONAG(X)$, for some set $Y_n\sub X\sm \Un_{k<n}V_{D_k}$.
Since the family $\cU_n$ covers $X\sm \Un_{k<n}V_{D_k}$, a family 
\[
\cU_n':=\sset{V^n_x\cup \Un_{k<n}V_{D_k}}{x\in Y_n}
\]
is an open cover of $X$.
Assume that $\cU_n'$ is the next move of \Alice{} in $\M(X)$.
Let $\cF_n'\sub\cU_n'$ be a finite set played by \Bob{} in $\M(X)$, according to $\sigma$, as a response to $\cU_n'$.
There is a finite set $D_n\sub Y_n$ such that 
\[
\cF'_n=\sset{V^n_x\cup \Un_{k<n}V_{D_k}}{x\in D_n}.
\]
We have $\Un\cF_n'=\Un_{k\leq n}V_{D_k}$.
Then in $\PONAG(X)$, \Bob{} plays the set $D_n$.

Since $\sigma$ is a winning strategy for \Bob{} in $\M(X)$, the family $\Un_{n\in\w}\cF'_n$ covers $X$.
Since $\Un\cF_n'=\Un_{k\leq n}V_{D_k}$ for each $n$, the family $\sset{V_{D_n}}{n\in\w}$ covers $X$.
It follows that the defined above strategy for \Bob{} in $\PONAG(X)$ is a winning strategy.
\epf

Recall that a space $X$ is Fr\'echet if for every $A\sub X$ and $p\in \overline{A} $ there exists a sequence in $A$ converging to $p$.

\bprp\label{prp:ponag-men}
Let $X$ be a Lindel\"of Fr\'echet space.
If \Bob{} has a winning strategy in $\PONAG(X)$, then \Bob{} has a winning strategy in $\M(X)$. 
\eprp

\bpf
Let $\sigma$ be a winning strategy for \Bob{} in $\PONAG(X)$.
We construct a winning strategy for \Bob{} in $\M(X)$ as follows.

{\bf Round 0:} Let ${\cU}_0$ be a family played by \Alice{} in $\M(X)$.
Since the space $X$ is Lindel\"of, we may assume that $\cU_0=\sset{U^0_k}{k\in\w}$ is an increasing cover of $X$.
We may also assume that there is a non-isolated point $y_0$ in $X$.
Indeed, because otherwise, the space is countable and we are done.
By the Fr\'echet property, there is a sequence $Y_0=\sset{y^0_k}{k\in\w}\sub X$ converging to $y_0$.
For simplicity we may assume that $y_0\in U_0^0$.
Let $V_{y_k^0}:=U_k^0$ for all $k$.
Assume that in $\PONAG(X)$ \Alice{} plays the family $\cU_0':=\sset{V_{y_k^0}}{y_k^0\in Y_0}$.
Let $D_0\sub Y_0$ be a set played by \Bob{} according to $\sigma$, in response to $\cU_0'$.
Since the set $D_0$ is a subsequence of a convergent sequence and it is closed in $X$, the set $D_0$ is finite.
There is a set $U_0\in\cU_0$ such that $U_0=V_{D_0}$ and \Bob{} plays the set $U_0$ in $\M(X)$ as a response to $\cU_0$.

{\bf Round 1:} Assume that ${\cU}_1$ is a family played by \Alice{} in $\M(X)$.
Since the space $X$ is Lindel\"of, we may assume that $\cU_1=\sset{U^1_k}{k\in\w}$ is an increasing cover of $X$ and also that $U_0\subseteq U_0^1$.
As before we may assume that there is a non-isolated point $y_1$ in $X\sm V_{D_0}$.
By the Fr\'echet property, there is a sequence $Y_1=\sset{y^1_k}{k\in\w}\sub X\setminus V_{D_0}$ converging to $y_1$.
Let $V_{y_k^1}:=U_k^1$ for all $k$.
Assume that in $\PONAG(X)$ \Alice{} plays the family $\cU_1':=\sset{V_{y^1_k}}{y^1_k\in Y_1}$.
Let $D_1\sub Y_1$ be a set picked by \Bob{} according to $\sigma$, in response to $\cU_1'$.
Since the set $D_1$ is a subsequence of a convergent sequence and it is closed in $X$, the set $D_1$ is finite.
There is a set $U_1\in\cU_1$ such that $U_1=V_{D_1}$ and \Bob{} plays the set $U_1$ in $\M(X)$ as a response to $\cU_1$.

Fix $n>1$ and assume that
\bi
\item families $\cU_k$ played by \Alice{} in $\M(X)$,
\item sets $U_k$ played by \Bob{} in $\M(X)$,
\item convergent sequences $Y_k:=\sset{y^k_i}{i\in\w}$
\item partial open neighborhood assignments $\sset{V_{y^k_i}}{y^k_i}\in Y_k$ played by \Alice{} in $\PONAG(X)$
\item sets $D_k\sub Y_k$ played by \Bob{} in $\PONAG(X)$ according to the strategy $\sigma$ in response to $\cU_k'$
\ei
such that $U_k=V_{D_k}$ have already been defined for all $0\leq k<n$ and $U_0\sub\dotsb\sub U_{n-1}$.

{\bf Round n:} 
Assume that $\cU_n$ is a family played by \Alice{} in $\M(X)$. Since the space $X$ is Lindel\"of, we may assume that $\cU_n=\sset{U^n_k}{k\in\w}$ is an increasing cover of $X$ and also that $U_{n-1}\subseteq U_0^n$.
We may assume that there is a non-isolated point $y_n$ in $X\sm V_{D_n}$.
By the Fr\'echet property, there is a sequence $Y_n=\sset{y^n_k}{k\in\w}\sub X\setminus V_{D_{n-1}}$ converging to $y_n$.
Let $V_{y_k^n}:=U_k^n$ for all $y_k^n$.
Assume that in $\PONAG(X)$, \Alice{} plays the family $\cU_n':=\sset{V_{y_k^n}}{y^k_n\in Y_n}$.
Let $D_n\sub Y_n$ be a set picked by \Bob{} according to $\sigma$, in response to $\cU_n'$.
Since the set $D_n$ is a subsequence of a convergent sequence and it is closed in $X$, the set $D_n$ is finite.
There is a set $U_n\in\cU_n$ such that $U_n=V_{D_n}$ and \Bob{} plays the set $U_n$ in $\M(X)$ as a response to $\cU_n$.

Since $\sigma$ is a winning strategy for $\Bob{}$ in $\PONAG(X)$, the family $\sset{V_{D_n}}{n\in\w}$ covers $X$. 
Since $U_n=V_{D_n}$ for all $n$, the family $\sset{U_n}{n\in\w}$ also covers $X$, and thus the defined above strategy for $\Bob{}$ in $\M(X)$ is a winning strategy.

\epf

\brem
The proof of Theorem~\ref{Alice_game_equivalence} and
Proposition~\ref{prp:ponag-men} remain true if \Bob{} has a winning strategy in a modified version of $\PONAG$, where in each round $n$, \Bob{} in response to an \Alice{}'s partial open neighborhood assignment $\sset{V_x}{x\in Y_n}$, plays just a closed subset $D_n\subseteq Y_n$ in $X$ (not necessarily closed and discrete in the ambient space $X$).
\erem

\bcor
Let $X$ be a Lindel\"of Fr\'echet space.
Then $\PONAG(X)$ and $\M(X)$ are equivalent games.
\ecor

 By a result of Telgarsky~\cite{Telg1}, if a space $X$ is perfectly normal and \Bob{} has a winning strategy in $\M(X)$, then the space $X$ is $\sigma$-compact. 
 Combination of this fact together with our results gives the following corollary.

\bcor
Let $X$ be a perfectly normal space. Then the following assertions are equivalent.
\be
\item \Bob{} has a winning strategy in $\M(X)$.
\item \Bob{} has a winning strategy in $\PONAG(X)$.
\item The space $X$ is $\sigma$-compact.
\ee
\ecor

\section{Comments and open problems}

\subsection*{Subspaces of the Sorgenfrey line and Lindel\"ofness in all finite powers}

By the result of Michael~\cite{michael}, assuming that \CH{} holds, there is an uncountable subspace of the Sorgenfrey line, whose all finite powers are Lindel\"of.
On the other hand, by the result of Todorcevic~\cite{tod}, if the open colorings axiom holds (a statement independent from ZFC), each uncountable subspace of the Sorgenfrey line has a non-Lindel\"of square.

A space $X$ is \emph{Hurewicz}~\cite{H}, if for each sequence $\eseq{\cU}$ of open covers of $X$, there are finite sets $\cF_0\sub\cU_0,\cF_1\sub\cU_1,\dotsc$ such that the family $\sset{\Un\cF_n}{n\in\w}$ is a \emph{$\gamma$-cover} of $X$, i.e., the sets $\sset{n}{x\in\Un\cF_n}$ are cofinite for all $x\in X$.
We have the following implications
\[
\sigma\text{-compactness}\longrightarrow \text{Hurewicz} \longrightarrow \text{Menger}
\]
and none of these implications is reversible, even in the class of subspaces of the real line~\cite{coc2, bst, chp, sfh}.

Let $\PN$ be the powerset of $\w$.
We identify each set from $\PN$ with its characteristic function, an element of the Cantor cube $\Cantor$.
Let $\roth$ be the family of all infinite subset of $\w$ and $\Fin$ be the family of all finite subsets of $\w$.
Identifying each set $x\in\roth$ with the increasing enumeration of its elements, we have that $x\in\NN$.
For $x,y\in \roth$, we write $x\les y$ if the set $\sset{n}{x(n)\leq y(n)}$ is cofinite.
A set $A\sub\roth$ is \emph{unbounded}, if for each $y\in \roth$, there is $x\in A$ such that $x\not\les y$.
Let $\fb$ be the minimal cardinality of an unbounded set in $\roth$.
We have $\w_1\leq\fb\leq\fc$ and it is independent from ZFC that $\fb=\w_1<\fc$.
It turns out that the above mentioned Michael result can be obtained using much milder hypothesis than \CH{}.

\bexm
Assume that $\fb=\w_1$.
There is an uncountable subspace of the Sorgenfrey line whose all finite powers are Hurewicz (in particular, they are  Lindel\"of).

The set $\PN$ can be treated as a subset of the real line using standard homeomorphism between $\Cantor$ and the standard middle-third Cantor set $C$ in the real line.
Using this standard homeomorphism and above identifications of elements in $\PN$ and $\Cantor$, we see
that a counterpart in $\PN$ of the set of all left-endpoints of intervals which are engaged in a construction of this standard middle-third Cantor set $C$ is the set $\Fin$.
Thus, we may assume that $\PN=C$.
Let $\PN_\text{M}$ be the set $\PN$ with a topology, where open neighborhoods of points from $\Fin$ are inherited from the standard topology on $\bbR$ and the points from $\roth$ are isolated.

Let $X=\sset{x_\alpha}{\alpha<\fb}\sub \roth$ be an unbounded set in $\roth$ such that $x_\alpha\les x_\beta$ for all $\alpha<\beta<\fb$ (such a set exists in ZFC).
Let $(X\cup\Fin)_\text{M}$ be the set $X\cup\Fin$ with a topology inherited from $\PN_\text{M}$.
Now let $(X\cup\Fin)_\text{S}$ be the set $X\cup\Fin$ with a topology inherited from the Sorgenfrey line.
Observe that in such a case basic open sets for the points from $\Fin$ are exactly the same, when consider $\Fin$ as a subspace of $(X\cup\Fin)_\text{M}$ and $(X\cup\Fin)_\text{S}$.
Thus, the identity map from $(X\cup\Fin)_\text{M}$ onto $(X\cup\Fin)_\text{S}$ is continuous.

By the result of Tsaban and the third named author~\cite{pMGen}, since $\fb=\w_1$, all finite powers of $(X\cup\Fin)_\text{M}$ are Hurewicz.
Since $(X\cup\Fin)_\text{S}$ is a continuous image of $(X\cup\Fin)_\text{M}$ and the Hurewicz property is preserved by continuous mappings, all finite powers of $(X\cup\Fin)_\text{S}$ are Hurewicz (and thus Lindel\"of).
\eexm

\subsection{Open problems}

We finish the paper with the following problems emerging from our investigations.
A wider class than $\sigma$-compact spaces is the class of Hurewicz spaces. 
In the light of Theorem~\ref{Th:tot_par_impl_Men1} the following problem arises.

\bprb
Let $X$ be a totally paracompact space with a dense Hurewicz subspace.
Is the space $X$ Menger?
What if this dense subset is a Lindel\"of space of size smaller than $\fb$? 
What if the Martin Axiom holds plus the negation of \CH{}?
\eprb

\bprb
Are the Menger game and the partial open neighborhood assignment game equivalent in the class of Lindel\"of spaces?
\eprb

\bprb
Let $X$ be a Menger space and assume that \Bob{} has a winning strategy in $\PONAG(X)$.
Does it follow that \Bob{} has a winning strategy in $\M(X)$?
\eprb

We would like also to know whether there is some relation between the partial open neighborhood assignment game and total paracompactness.

\bprb
Let $X$ be a regular space.
\be
\item Does \Alice{} has no winning strategy in $\PONAG(X)$ imply that $X$ is totally paracompact?
\item Does \Bob{} has a winning strategy in $\PONAG(X)$ imply that $X$ is totally paracompact?
\ee
\eprb

\section*{Acknowledgments}
The authors thank the National Group for Algebric and Geometric Structures, and their Applications (GNSAGA-INdAM) for its invaluable support.
Parts of the work
reported here were carried out during two visits of the third named author at the University of Messina.
These visits were partially supported by the University of Messina. The third named author
thanks the other authors for their kind hospitality.


\begin{thebibliography}{99}
	
	\bibitem{Aurichi} L. Aurichi, \emph{$D$-spaces, topological games, and selection principles}, Topology Proceedings \textbf{36} (2010), 107--122.

\bibitem{ABG} L. Aurichi, M. Bonanzinga, D. Giacopello, \emph{On some topological games involving networks}, Topology and its Applications 351 (2024), Art. ID 108936.
    
	\bibitem{BB}  Z. Balogh, H. Bennett, \emph{Total paracompactness of real GO-spaces}, Proceedings of the American Mathematical Society 101(1987), 753–-760.

\Pa{bst}{T. Bartoszy\'nski, S. Shelah, B. Tsaban}{Additivity Properties of Topological Diagonalizations}{Journal of Symbolic Logic}{68}{2003}{1254}{1260}

\bibitem{chp} J. Chaber, R. Pol, \emph{A remark on Fremlin–Miller theorem concerning
the Menger property and Michael concentrated sets}, available at: https://arxiv.org/pdf/2305.10797
	
\bibitem{B}
C. Bandy, \emph{Property Q}, International Journal of Mathematics and Mathematical Sciences, 14 (1991), 315‑-318.
	
\bibitem{Curtis} D. Curtis, \emph{Total and absolute paracompactness}, Fundamenta Mathematicae 77 (1973), 277-–283.

\bibitem{D} E. van Douwen, W. Pfeﬀer, \emph{Some properties of the Sorgenfrey line and related spaces}, Pacific Journal of Mathematics 81 (1979), 371–-377.

\bibitem{engelking} R. Engelking, \emph{General Topology}, Sigma Series in Pure Mathematics, Heldermann Verlag (1989).
	
	\bibitem{Ford} R. Ford, \emph{Basis properties in dimension theory}, Doctoral Dissertation, Auburn University,
	Auburn, AL, 1963.
	
\bibitem{H} W. Hurewicz, \"Uber eine Verallgemeinerung des Borelischen Theorems, Mat. Z. 24 (1926) 401--421.

\bibitem{coc2} W. Just, A. Miller, M. Scheepers, P. Szeptycki,
\emph{The combinatorics of open covers II},
Topology and its Applications \textbf{73} (1996), 241--266.

\bibitem{lelek}A. Lelek, \emph{Some cover properties of spaces}, Fundamenta Mathematicae 64(1964), 209–-218.

\bibitem{michael}E. Michael, \emph{Paracompactness and the Lindel\"of property in finite and countable Cartesian products}, Compositio Mathematica 23 (1971), 119–-214.

\bibitem{Sakai} M. Sakai, \emph{Menger subsets of the Sorgenfrey line}, Proceedings of the American Mathematical Society 137 (2009), 3129–-3138.

\bibitem{SteenSeebach} L.A. Steen and J.A. Seebach, \emph{Couterexamples in Topology}, Springer-Verlag,
Berlin - Heidelberg - New York, 1978.

\bibitem{SzTsGame} P. Szewczak, B. Tsaban, \emph{Conceptual proofs of the Menger and Rothberger games}, Topology and its Applications 272 (2020), 107048.

\bibitem{pMGen} P. Szewczak, B. Tsaban, \emph{Products of general Menger spaces}, Topology and its Applications 255 (2019), 41--55.

\bibitem{telg}R. Telg\'{a}rsky, \emph{Closure‑preserving covers}, Fundamenta Mathematicae 85(1974), 165-‑175.

\bibitem{Telg1} R. Telg\'arsky, \textit{On games of Tops{\o}e}, 	Mathematica Scandinavica 54 (1984).

\bibitem{tod} S. Todorcevic, \emph{Partition problems in topology}, Contemporary Math. 84. Amer. Math. Soc. Providence,
RI, (1989). xii+116 pp.

\Pa{sfh}{B. Tsaban, L. Zdomskyy}{Scales, fields, and a problem of Hurewicz}{Journal of the European Mathematical Society}{10}{2008}{837}{866}

	
	
\end{thebibliography}
\end{document}